# New Solution Principles of Multi-criteria Problems Based on Comparison Standards

## V.O. Groppen

*(Vladikavkaz, Russia)*

Proposed is a new formal approach for solution of extreme multi-criteria problems transforming them into single-criterion mathematical models, without any additional information. Transforming rules are based on comparison standards and intervals between these standards and goal function values, corresponding to problem solution. Pareto-optimality of these solutions is proved by a number of theorems. Examples illustrating efficiency and universality of this approach are presented.

### 1. Introduction

As a rule search for multi-criteria problems solution is based on Pareto optima principle [8]: vector of variables is optimal if improvement of values of any subset of goal functions is always accompanied by a change of another goal functions subset for the worse. One of the most important imperfections of such a definition is its' comparatively low selectivity: power of set of optimal solutions is often commensurable with a set of all feasible plans transforming it into a sampling problem [3]. To avoid this, a multi-criteria problem is usually transformed into a problem with a single goal function with the help of certain additional information, which is not available in the original problem statement. Examples of such approach are lexicographic goals ordering and creation of super-criterion generated as a sum of weighted goal functions of original problem [4, 6, 7]. Contrary to them transformation of any multi-criteria problem into a problem with a single goal function that does not need any additional condition with the use of comparison standards is proposed [4]. These standards are connected with combinations of the best (comparison standard "a") and the worst (comparison standard "b") values of goal functions regardless of the existence of corresponding feasible vectors of variables. It is possible to say that point "b" corresponding to any multi-criteria problem personifies the absolute evil, whereas point "a" can be personified for this problem as the absolute good.

Presented below analysis is based on search for such a feasible vector of variables which corresponds to a point in criteria space with extreme intervals to comparison standards.

### 2. Transformation of multi-objective optimization problems

Below the following symbols and definitions are used:

$\vec{X}$ - vector of variables;

$X_k$ – set of k - th variable values,

$F_i(\vec{X})$ – i-th criterion (i = 1, 2, ..., n);

I – set of indices used for criteria enumeration ($|I| = n$);

$\varphi_j(\vec{X})$ – j-th restriction of an extreme problem,

$K_i$ – numerical estimate of the best value for i-th criterion;

$W_i$ - numerical estimate of the worst value for i-th criterion;

G(X, U) – weighted directed graph with set of vertices X and set of arcs U;
r(i, j) – capacity of arc (i, j) ∈ U;
A(G) – set of circuits in G(X, U);

$a_i$ - i-th circuit in set À(G);

$y_i$ – circulation ( a flow in a circuit with coinciding source and terminal vertices) value in
$a_i \in A(G)$.

Below circulation in graph G(X, U) is defined as a set of circulations in circuits of set A [1,2].

It is obvious that value of each i-th criterion is neither better than $K_i$, nor worse than $W_i$: $\forall i,\ K_i \succ F_i(\vec{X}) \succ W_i$ .

Situation when values of $K_i$ and $W_i$ coincide is not analyzed because it is a case of $F_i(\vec{X})$ constant.

Formal statement of multi-criteria problem is given below as follows:

(1) $\begin{cases} \forall i: F_i(\vec{X}) \to \max\ (\min); \\ \forall j: \varphi_j(\vec{X}) \le b_j; \\ \forall k: x_k \in X_k;\ \vec{X} = \{x_1,\ x_2,\ .....,\ x_m.\}. \end{cases}$

It permits us to determine values of $K_i$ (i ∈ I) solving n systems:

(2) $\begin{cases} K_i = F_i(\vec{X}) \to \max\ (\min); \\ \forall j: \varphi_j(\vec{X}) \le b_j; \\ \forall k: x_k \in X_k;\ \vec{X} = \{x_1,\ x_2,\ .....,\ x_m.\}. \end{cases}$

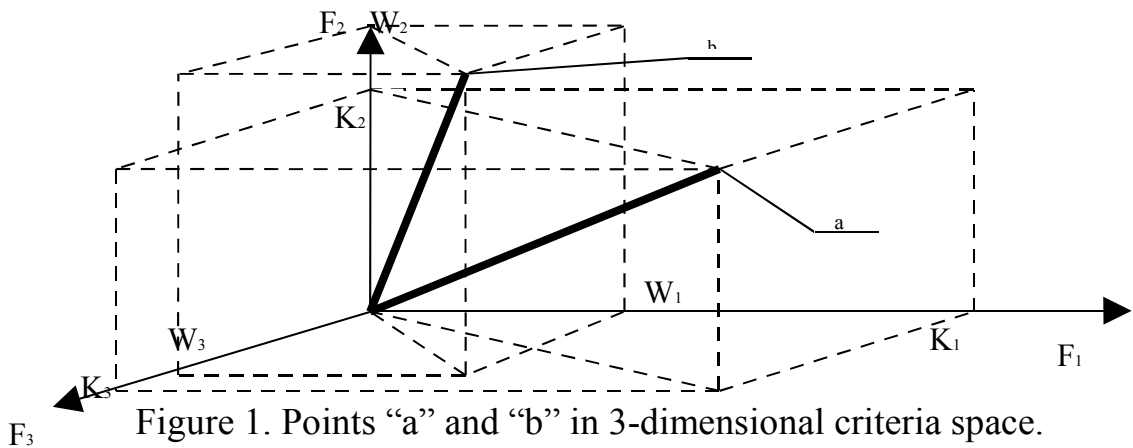

Figure 1. Points "a" and "b" in 3-dimensional criteria space.

Solution of (2) for each $K_i$, ($i \in I$), results in vector $\vec{X}$:

(3) $\quad \vec{K} = \{K_1, K_2, ....., K_n\}$,

corresponding to point "a" in the n-dimensional criteria space (Figure 1).

Using similar technology with the inverse goal for each criterion it is possible to determine in the same space point "b" (Figure 1) with corresponding vector $\vec{W}$ = { W$_1$, W$_2$, .., Wn}, reflecting combination of the worst values of criteria $F_i(\vec{X})$, $i \in I$, for system (1). Evidently either one or both of them can correspond to a vector of variables which is not feasible.

Approaches to transformation of system (1) to the solution of problem with a single goal function:
1. Searched is system (1) feasible vector of variables corresponding in the criteria space to point "m" with minimal interval Δ(a, m) between "m" and "a".
2. Searched is system (1) feasible vector of variables corresponding to a point "n" with maximal interval θ(a, m) between "n" and "b" in the criteria space.
3. Searched is system (1) feasible vector of variables corresponding to a point "d" with minimal ratio value Δ(a, m)/θ(a, m) in the criteria space.

## 2.1. Multi-criteria problems transforming in the case of homogeneous criteria

If all the criteria are homogeneous, then interval between any two points in criteria space is determined as a distance in Euclidean space [7], the last being used for system (1) transformation into a single-criterion problem. Thus searching system (1) feasible vector of variables corresponding in the criteria space to the point "m" with minimal interval Δ(a, m) between "m" and "a" results in minimizing of the following function:

(4) $\quad \Delta = \sqrt{\sum_i [K_i - F_i(\vec{X})]^2}$.

It permits to substitute (1) by the following single-criterion system:

(5) $\quad \begin{cases} \Delta = \sqrt{\sum_{i=1}^{n} [K_i - F_i(\vec{X})]^2} \to \min; \\ \forall j : \varphi_j(\vec{X}) \leq b_j; \\ \forall k : x_k \in X_k; \quad \vec{X} = \{x_1, x_2, ....., x_m.\}. \end{cases}$

True is the following theorem, reflecting connection between systems (1) and (5):

**Theorem 1.** System (5) optimal vector of variables is simultaneously Pareto-optimal system (1) solution.

Proof of theorem 1 is presented in the Appendix.

Similarly searching system (1) feasible vector of variables, reflecting in the criteria space by search of the point "n" with maximal interval θ(n, b) between "n" and "b", corresponds to system (1) substituted by the following one:

$$(6) \begin{cases} \theta = \sqrt{\sum_i [W_i - F_i(\vec{X})]^2} \to \max; \\ \forall j : \varphi_j(\vec{X}) \leq b_j; \\ \forall k : x_k \in X_k; \quad \vec{X} = \{x_1, x_2, ....., x_m\}. \end{cases}$$

For systems (1) and (6) the theorem, similar to theorem 1 above is true:

**Theorem 2.** System (6) optimal vector of variables is Pareto-optimal solution for system (1).

Proof of theorem 2 is presented in the Appendix.

Searching system (1) feasible vector of variables corresponding in the criteria space to "d" point with minimal interval Δ(a, d) between "d" and "a" and, simultaneously, with maximal interval θ(b, d) between "b" and "d" results in combination of systems (5) and (6):

$$(7) \begin{cases} \Delta = \sqrt{\sum_i [K_i - F_i(\vec{X})]^2} \to \min; \\ \theta = \sqrt{\sum_i [W_i - F_i(\vec{X})]^2} \to \max; \\ \forall j : \varphi_j(\vec{X}) \leq b_j; \\ \forall k : x_k \in X_k; \quad \vec{X} = \{x_1, x_2, ....., x_m\}, \end{cases}$$

Due to the fact that (7) is a two-criteria system, it is associated with only two lexicographic goals Δ and θ orderings: $\pi_1 = \{\Delta, \theta\}$, $\pi_2 = \{\theta, \Delta\}$. An important link between systems (1) and (7) solutions is given by the following theorem:

**Theorem 3.** System (7) Pareto-optimal vector of variables found for any lexicographic criteria order is simultaneously Pareto-optimal for system (1).

Proof of theorem 3 is presented in the Appendix.

To arrive at an optimal solution of system (7), the latter can be transformed into the following single-criterion problem:

(8)
$$\begin{cases} \gamma = \sqrt{\dfrac{\sum_i [K_i - F_i(\vec{X})]^2}{\sum_i [W_i - F_i(\vec{X})]^2}} \to \min; \\ \forall j : \varphi_j(\vec{X}) \le b_j; \\ \forall k : x_k \in X_k; \quad \vec{X} = \{x_1, x_2, \ldots, x_m\}, \end{cases}$$

due to the following theorem:
**Theorem 4.** System (8) optimal vector of variables is simultaneously Pareto-optimal system (7) solution.
Proof of theorem 4 is given in the Appendix.
  Thus, due to theorems 1 – 4 search of Pareto-optimal multi-criteria problem solution can be replaced by that of a single-criterion problem. Example illustrating this approach is presented below.
**Example 1.** In a directed, weighted and strongly connected graph G(X, U) (Figure 2) with capacity of each arc equal to 1, maximal circulation in each simple circuit is searched [1 – 3].

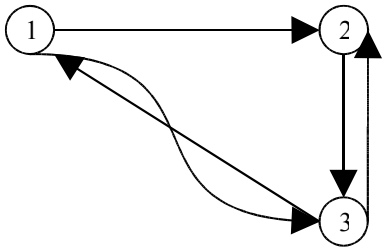

Figure 2. Graph G(X,U)

Circuits set A(G) consists of three simple circuits [7]:
$a_1 = \{x_1, (1,3), x_3, (3,1), x_1\}$;
$a_2 = \{x_2, (2,3), x_3, (3,2), x_2\}$;
$a_3 = \{x_1, (1,2), x_2, (2,3), x_3, (3,1), x_1\}$,
thus permitting the following problem formal statement:

(9)
$$\begin{cases} \forall i : y_i \to \max; \\ y_1 + y_3 \le 1; \\ y_2 + y_3 \le 1; \\ \forall i : 0 \le y_i \le 1. \end{cases}$$

Keeping in mind (9), system (8) is transformed as follows:

(10)
$$\begin{cases} \gamma = \sqrt{\dfrac{\sum_{i=1}^{n}(1-y_i)^2}{\sum_{i=1}^{n}y_i^2}} \to \min; \\ y_1 + y_3 \le 1; \\ y_2 + y_3 \le 1; \\ \forall i: 0 \le y_i \le 1. \end{cases}$$

Problem (10) optimal solution $y_1 = y_2 = 1,\ y_3 = 0,$ coincides with one of system (9) Pareto-optimal solutions. It is possible to show that substituting system (10) goal function by minimized criteria (4), and keeping in mind that $\forall i, K_i = 1$, results in the other Pareto-optimal system (9) solution: $\forall i,\ y_i = 0.5$.

It can be easily seen that elimination of square root in (5) – (8) goal functions results in the same optimal vectors of variables, as in an initial problem statements. Possibility to use similar approach for the case of heterogeneous criteria is proved below.

## 2.2. Case of heterogeneous criteria

As it was mentioned, the approach described above operates with Euclidean distances in the space of homogeneous criteria. Heterogeneous criteria cause heterogeneous space of these criteria, thus restricting usage of developed approach.

Normalization of criteria permits us to avoid this restriction, substituting (1) by system:

(11)
$$\begin{cases} \forall i: f_i(\vec{X}) = \dfrac{F_i(\vec{X}) - \min(K_i; W_i)}{\max(K_i; W_i) - \min(K_i; W_i)} \to \max\ (\min); \\ \forall j: \varphi_j(\vec{X}) \le b_j; \\ \forall k: x_k \in X_k;\quad \vec{X} = \{x_1,\ x_2,\ \ldots,\ x_m.\}, \end{cases}$$

with dimensionless criteria values and with spread of these values 0 – 1 for each criterion. Connecting the above approach with multi-criteria problem with heterogeneous criteria solution, the following holds:

**Theorem 5.** System (11) Pareto-optimal vector of variables is simultaneously system (1) Pareto-optimal solution.

Proof of theorem 5 is given in the Appendix.

**Example 2.** In a directed, weighted and strongly connected graph G(X, U) (Figure 2) with the capacity of each arc equal to 1, we search such a circulation, that:
- the sum of circulations in simple circuits is maximized;
- the total cost of circulations in simple circuits is minimized;
- the cost of a single unit of circulation in i-th circuit is equal to (i – 1);

- the total circulation in $(i,j) \in U$ arc belonging to several simple circuits is equal to
  $r(i, j)=1$.

Formal statement of this problem looks as follows:

(12)
$$\begin{cases} F_1 = \sum_i y_i \to \max; \\ F_2 = \sum_i (i-1) y_i \to \min; \\ y_1 + y_3 = 1; \\ y_2 + y_3 = 1; \\ \forall i : 0 \leq y_i \leq 1. \end{cases}$$

Keeping in mind that for this problem $K_1 = 2$, $W_1 = 1$, $K_2 = 1$, $W_2 = 2$, system (12) transformed and described as (11) looks like:

(13)
$$\begin{cases} \gamma = \sqrt{\dfrac{[1-F_1]^2 + F_2^2}{F_1^2 + [1-F_2]^2}} = \sqrt{\dfrac{[1-\sum_i y_i]^2 + [\sum (i-1)y_i]^2}{[\sum_i y_i]^2 + [1-\sum (i-1)y_i]^2}} \to \min; \\ y_1 + y_3 = 1; \\ y_2 + y_3 = 1; \\ \forall i : 0 \leq y_i \leq 1. \end{cases}$$

It is easy to see that system (13) optimal vector of variables: $y_1 = y_2 = 1$, $y_3 = 0$, is, simultaneously, Pareto-optimal solution for system (12) with criteria values $F_1 = 2$, $F_2 = 1$.

## 3. Conclusion

Presented above approach differs from the usually used for multi-criteria problems solution by its transformation into a single-criterion form with the following features:
- it needs no additional information, being absent in the original problem statement;
- it guarantees Pareto-optimality of initial problem solution.

## 5. Appendix

**Proof of Theorem 1.**

Let us denote system (5) optimal vector of variables as $\vec{X}_{opt}$, and using the rule of contraries assuming that the theorem is wrong, i.e. there is a Pareto-optimal vector of variables $\vec{X}_1$ dividing the set of criteria **F** into two subsets $\mathbf{F}_1$ and $\mathbf{F}_2$ with the following properties:

(14)    $\forall F_i \in \mathbf{F}_1 : F_i(\vec{X}_1) = F_i(\vec{X}_{opt})$;

(15)    $\forall F_i \in \mathbf{F}_2 : [F_i(\vec{X}_1) - F_i(\vec{X}_{opt})]^2$.

Using (14) and (15) difference of squares for system (5) goal function is determined as follows:

(16)    $\Delta^2(\vec{X}_{opt}) - \Delta^2(\vec{X}_1) = \sum_{F_i \in F_1}[K_i - F_i(\vec{X}_{opt})]^2 - \sum_{F_i \in F_1}[K_i - F_i(\vec{X}_1)]^2 + \sum_{F_i \in F_2}[K_i - F_i(\vec{X}_{opt})]^2 - \sum_{F_i \in F_2}[K_i - F_i(\vec{X}_1)]^2$.

As, due to condition (14), the difference of first two sums of the second member of equality (16), is equal to zero, this equality can be replaced by the following:

(17)    $\Delta^2(\vec{X}_{opt}) - \Delta^2(\vec{X}_1) = \sum_{F_i \in F_2}[K_i - F_i(\vec{X}_{opt})]^2 - \sum_{F_i \in F_2}[K_i - F_i(\vec{X}_1)]^2$

After removal of brackets in the second member of equality (17), the last can be presented as sum of products:

(18)    $\Delta^2(\vec{X}_{opt}) - \Delta^2(\vec{X}_1) = \sum_{F_i \in F_2}[F_i(\vec{X}_1) - F_i(\vec{X}_{opt})][F_i(\vec{X}_1) + F_i(\vec{X}_{opt}) - 2K_i]$

Due to the definitions and assumptions made, there are only two possibilities for $K_i, F_i(\vec{X}_1)$ and $F_i(\vec{X}_{opt})$ correlation :

1. If $F_i(\vec{X}) \to \min$, than the following holds : $K_i \leq F_i(\vec{X}_{opt}) \leq F_i(\vec{X}_1)$.

2. If $F_i(\vec{X}) \to \max$, then true is the following inequality $K_i \geq F_i(\vec{X}_{opt}) \geq F_i(\vec{X}_1)$.

It is easy to see that each case above results in nonnegative second member of equality (18):

$$\Delta^2(\vec{X}_{opt}) - \Delta^2(\vec{X}_1) = [\Delta(\vec{X}_{opt}) - \Delta(\vec{X}_1)][\Delta(\vec{X}_{opt}) + \Delta(\vec{X}_1)] \geq 0.$$

As value $\Delta$ is nonnegative for all possible combinations of values of vector of variables components, true is inequality:

(19) $\quad \Delta(\vec{X}_{opt}) - \Delta(\vec{X}_1) \geq 0.$

In other words (19) proves that vector $\vec{X}_{opt}$ is not optimal. This result is at variance with the accepted above assumption, which was to be proved.

**Proof of Theorem 2.**

Let us denote system (9) optimal vector of variables as $\vec{X}_{opt}$, and using the rule of contraries assume that the theorem is wrong, i.e. there is such a Pareto-optimal vector of variables $\vec{X}_1$, that:

a) (20) $\quad \theta(\vec{X}_{opt}) - \theta(\vec{X}_1) \geq 0;$

b) set of criteria **F** can be divided into two subsets **F**¹ and **F**² with components fulfilling conditions (14) – (15).

Using a) and b) it's possible to determine squares difference for system (9) goal function as follows:

(21) $\quad \theta^2(\vec{X}_{opt}) - \theta^2(\vec{X}_1) = \sum_{F_i \in F_1}[W_i - F_i(\vec{X}_{opt})]^2 - \sum_{F_i \in F_1}[W_i - F_i(\vec{X}_1)]^2 + \sum_{F_i \in F_2}[W_i - F_i(\vec{X}_{opt})]^2 - \sum_{F_i \in F_2}[W_i - F_i(\vec{X}_1)]^2$

Due to (14) equality (21) is transformed into:

(22) $\quad \theta^2(\vec{X}_{opt}) - \theta^2(\vec{X}_1) = \sum_{F_i \in F_2}[F_i(\vec{X}_{opt}) - F_i(\vec{X}_1)][F_i(\vec{X}_1) + F_i(\vec{X}_{opt}) - 2W_i].$

Thanks to the definitions and assumptions made, there are only two possibilities for correlation of $W_i, F_i(\vec{X}_1)$ and $F_i(\vec{X}_{opt})$:

- If $F_i(\vec{X}) \to \min$, then the following inequalities are true :

(23) $\quad W_i \leq F_i(\vec{X}_{opt}) \leq F_i(\vec{X}_1)$.

- If $F_i(\vec{X}) \to \max$, then opposite inequalities are true :

(24) $$W_i \geq F_i(\vec{X}_{opt}) \geq F_i(\vec{X}_1).$$

Keeping in mind (23) and (24) it is easy to see that the second member of equality (22) is not positive, thus resulting in inequality:

(25) $$[\theta(\vec{X}_{opt}) + \theta(\vec{X}_1)][\theta(\vec{X}_{opt}) - \theta(\vec{X}_1)] \leq 0.$$

As expression in the first square brackets is nonnegative, true is inequality:

(26) $$\theta(\vec{X}_{opt}) - \theta(\vec{X}_1) \leq 0.$$

As inequality (26) conflicts with (20), assumption that the theorem is wrong fails. That was what we aimed to prove.

**Proof of Theorem 3.**

Below consequently both criteria permutations $\pi_1$ and $\pi_2$ are analyzed.

1. $\pi_1 = \{\Delta, \theta\}$. Eliminating criterion $\theta$ in system (7), we transform this system into (5), keeping in mind that due to theorem 1 optimal vector of variables of the latter $\vec{X}_{opt1}$ is simultaneously Pareto-optimal solution for system (1). Inserting equality:

(27) $$\Delta(\vec{X}_{opt1}) = \sqrt{\sum_i [K_i - F_i(\vec{X})]^2},$$

into system (6) the latter is transformed into the following system:

(28) $$\begin{cases} \theta = \sqrt{\sum_i [W_i - F_i(\vec{X})]^2} \to \max; \\ \sqrt{\sum_i [K_i - F_i(\vec{X})]^2} = \sqrt{\sum_i [K_i - F_i(\vec{X}_{opt1})]^2}; \\ \forall j: \varphi_j(\vec{X}) \leq b_j; \\ \forall k: x_k \in X_k; \quad \vec{X} = \{x_1, x_2, \ldots, x_m\}, \end{cases}$$

with optimal solution determined by vector of variables $\vec{X}_{opt2}$ and coinciding with system (7) optimal solution for $\pi_1$ criteria order. As $\Delta(\vec{X}_{opt2}) = \Delta(\vec{X}_{opt1})$, and vector of variables $\vec{X}_{opt2}$ satisfies system (5) restrictions, this vector also presents optimal solution of system (5), thus, due to theorem 1, Pareto-optimality for system (1) is obtained.

2. $\pi_2 = \{\theta, \Delta\}$. Eliminating criterion $\Delta$ in system (7), we transform this system into (6), keeping in mind that due to theorem 2 optimal vector of variables of the latter $\vec{X}_{opt3}$ is simultaneously Pareto-optimal solution for system (1). Inserting equality:

(29) $$\sqrt{\sum_i [W_i - F_i(\vec{X})]^2} = \sqrt{\sum_i [W_i - F_i(\vec{X}_{opt3})]^2},$$

into system (5) the last is transformed into the following system:

$$(30) \quad \begin{cases} \Delta = \sqrt{\sum_{i=1}^{n}[K_i - F_i(\vec{X})]^2} \to \min; \\ \sqrt{\sum_{i=1}^{n}[W_i - F_i(\vec{X})]^2} = \sqrt{\sum_{i=1}^{n}[W_i - F_i(\vec{X}_{opt3})]^2}; \\ \forall j: \varphi_j(\vec{X}) \le b_j; \\ \forall k: x_k \in X_k; \quad \vec{X} = \{x_1, x_2, \ldots, x_m\}. \end{cases}$$

By definition system (30) optimal vector of variables $\vec{X}_{opt4}$ coincides with system (7) optimal solution for criteria order $\pi_2$. As $\theta(\vec{X}_{opt3}) = \theta(\vec{X}_{opt4})$, and vector of variables $\vec{X}_{opt4}$ satisfies system (6) restrictions, this vector also presents optimal solution of system (6), thus, due to theorem 2, we obtain Pareto-optimality for system (1).
That was what we aimed to prove.

**Proof of Theorem 4.**
Proof is based on the rule of contraries. Denoting system (8) optimal vector of variables as $\vec{X}_{opt1}$, we suppose that the theorem is wrong, i.e. there is such a system (7) Pareto-optimal vector of variables $\vec{X}_{opt2}$, that fulfilled is one of the following conditions:

$$(31) \quad \begin{cases} \Delta(\vec{X}_{opt1}) > \Delta(\vec{X}_{opt2}); \\ \theta(\vec{X}_{opt1}) = \theta(\vec{X}_{opt2}); \end{cases}$$

$$(32) \quad \begin{cases} \Delta(\vec{X}_{opt1}) = \Delta(\vec{X}_{opt2}); \\ \theta(\vec{X}_{opt1}) < \theta(\vec{X}_{opt2}); \end{cases}$$

$$(33) \quad \begin{cases} \Delta(\vec{X}_{opt1}) > \Delta(\vec{X}_{opt2}); \\ \theta(\vec{X}_{opt1}) < \theta(\vec{X}_{opt2}); \end{cases}$$

Condition (31) entry results in the following inequality:

$$(34) \quad \frac{\Delta(X_{opt1})}{\theta(X_{opt1})} > \frac{\Delta(X_{opt2})}{\theta(X_{opt2})}.$$

In other words (34) means that $\gamma(\vec{X}_{opt1}) > \gamma(\vec{X}_{opt2})$, i.e. vector $\vec{X}_{opt1}$ is not optimal with reference to system (8) thus being at variance with the initial condition. Similarly it is easy to bring each condition (32), (33) to inequality (34), each time breaking initial condition of vector $\vec{X}_{opt1}$ optimality. It proves falseness of assumption that the theorem is wrong, what we aimed to prove.

**Proof of Theorem 5.**

Let us denote system (11) Pareto-optimal vector of variables as $\vec{X}_1$, and suppose that the theorem is wrong, i.e. there are vector of variables $\vec{X}_2$ and subset of indices $I' \subseteq I$, $I' \neq \varnothing$, with the following properties:

(35)
$$\begin{cases} \forall i \in I': f_i(\vec{X}_1) \succ f_i(\vec{X}_2); \\ \forall i \in I': F_i(\vec{X}_1) \prec F_i(\vec{X}_2); \\ \forall i \notin I': F_i(\vec{X}_1) = F_i(\vec{X}_2); \\ \forall i \notin I': f_i(\vec{X}_1) = f_i(\vec{X}_2); \end{cases}$$

It is easy to see that inserting the following expression for $f_i(\vec{X}_1)$ in (35):

$$\forall i: f_i(\vec{X}) = \frac{F_i(\vec{X}) - \min(K_i; W_i)}{\max(K_i; W_i) - \min(K_i; W_i)},$$

results in system (36):

(36)
$$\begin{cases} \forall i \in I': F_i(\vec{X}_1) \succ F_i(\vec{X}_2); \\ \forall i \in I': F_i(\vec{X}_1) \prec F_i(\vec{X}_2); \\ \forall i \notin I': F_i(\vec{X}_1) = F_i(\vec{X}_2); \\ \forall i \notin I': f_i(\vec{X}_1) = f_i(\vec{X}_2), \end{cases}$$

being simultaneous only for I' = ∅. The latter results in falseness of the above assumption, what we aimed to prove.